\documentclass[twoside, 12pt]{article}

\usepackage{filecontents,amsfonts,amsmath,amssymb,blkarray,multirow,graphics,amsthm,varioref,mathrsfs,enumerate}
\usepackage[noadjust]{cite}
\usepackage[colorlinks=true,allcolors=blue]{hyperref}
\usepackage[nameinlink,capitalize,noabbrev]{cleveref}
\usepackage{pgf,tikz}
\usetikzlibrary{arrows}
\crefname{equation}{}{}

\allowdisplaybreaks
\definecolor{wrwrwr}{rgb}{0.3803921568627451,0.3803921568627451,0.3803921568627451}
\definecolor{rvwvcq}{rgb}{0.08235294117647059,0.396078431372549,0.7529411764705882}
\newtheorem{theorem}{Theorem}[section]
\newtheorem{remark}[theorem]{Remark}
\newtheorem{example}[theorem]{Example}
\newtheorem{lemma}[theorem]{Lemma}
\newtheorem{corollary}[theorem]{Corollary}
\newtheorem{definition}[theorem]{Definition}
\newtheorem{proposition}[theorem]{Proposition}

\newtheorem{conjecture}{Conjecture}

\newcommand{\bea}{\begin{eqnarray}}
\newcommand{\eea}{\end{eqnarray}}

\newcommand{\comment}[1]{}

\def \bpm{\begin{pmatrix}}
\def \epm{\end{pmatrix}}
\def \bd{\begin{definition}}
\def \ed{\end{definition}}
\def \bcc{\begin{conjecture}}
\def \ecc{\end{conjecture}}
\def \bt{\begin{theorem}}
\def \et{\end{theorem}}
\def \bl{\begin{lemma}}
\def \el{\end{lemma}}
\def \bc{\begin{corollary}}
\def \ec{\end{corollary}}
\def\be#1\ee{\begin{align}#1\end{align}}
\def\beq #1\eeq {\begin{align*}#1\end{align*}}
\def \ben{\begin{enumerate}}
\def \een{\end{enumerate}}
\def \ba{\begin{array}}
\def \ea{\end{array}}
\def \bp{\begin{proposition}}
\def \ep{\end{proposition}}
\def \bx{\begin{example}}
\def \ex{\end{example}}
\def \br{\begin{remark}}
\def \er{\end{remark}}
\def \bdsc{\begin{description}}
\def \edsc{\end{description}}
\def\pf{{\noindent \it \bf Proof. }}
\def \qed {\hfill \vrule height6pt width6pt depth0pt}
\def\hs{\hspace{.3cm}}
\def\vs{\vskip .3cm}
\def\ds{\displaystyle}
\def\1{1\!\!1}
\def\J{\mathbb{J}}

\def\r{\rho}
\def\D{\mathcal{D}}
\def\M{\mathbb{M}}

\title{On Cartesian product of matrices}

\author{Deepak Sarma \\Department of Mathematical Sciences, \\ Tezpur University, Tezpur-784028, India.\\
Email addresses: \url{deepaks@tezu.ernet.in}}

\date{}

\begin{document}
\maketitle

\vs

\begin{abstract}
Recently, Bapat and Kurata [\textit{Linear Algebra Appl.}, 562(2019), 135-153] defined the Cartesian product of two square matrices $A$ and $B$ as $A\oslash B=A\otimes \J+\J\otimes B$, where $\J$ is the all one matrix of appropriate order and $\otimes$ is the Kronecker product. In this article, we find the expression for the trace of the Cartesian product of any finite number of square matrices in terms of traces of the individual matrices. Also, we establish some identities involving the Cartesian product of matrices.
Finally, we apply the Cartesian product to study some graph-theoretic properties.\\
\vs

\end{abstract}

\maketitle
\bigskip

\section{Introduction and terminology}
By $\mathbb{M}_{m,n}$, we denote the class of all matrices of size $m\times n.$ Also, by $\mathbb{M}_{n}$, we denote the class of all square matrices of order $n.$ For $M\in \mathbb{M}_{n}$ we write  $m_{ij}$ or $M_{ij}$ to denote the $ij-$th element of $M.$  By $\J $ and $\1$, we mean the matrix of all one's and vector of all one's, respectively
 of suitable order. Similarly $\large0$ denotes the zero matrix or the vector.  We will mention their order wherever its necessary.
Throughout this article, we denote the sum of all entries of a matrix $A$ by $S_A$ and the sum of the entries of $i-$th row of $A$ by $A_i.$ The inertia of a square matrix $M$ with real eigenvalues is the triplet $(n_+(M),n_0(M), n_-(M)),$ where $n_+(M)$ and $n_-(M)$ denote the number of positive and negative eigenvalues of $M$, respectively, and $n_0(M)$ is the algebraic multiplicity of 0 as an eigenvalue of $M.$ 

The \textit{Kronecker product} of two matrices $A$ and $B$ of sizes $m\times n$ and $p\times q,$ respectively, denoted by $A\otimes B$ is defined to be the $mp\times nq$ block matrix
\[A\otimes B= \bpm 
&a_{1,1}B &a_{1,2}B &\cdots &a_{1,n}B\\
&a_{2,1}B &a_{2,2}B &\cdots &a_{2,n}B\\
&\vdots &\vdots &\ddots &\vdots \\
& a_{m,1}B &a_{m,2}B &\cdots &a_{m,n}B
\epm .\] 

The \textit{Hadamard product} of two matrices $A$ and $B$ of the same size, denoted by $A\circ B$ is defined to be the entrywise product $A\circ B= [a_{i,j}b_{i,j}].$

Bapat and Kurata \cite{bk19} defined the \textit{Cartesian product} of two square matrices $A\in \M_m$ and $B\in \M_n$ as $A\oslash B=A\otimes \J_n+\J_m\otimes B.$ The authors proved the Cartesian product to be associative. We use $A^{[k]}$ to mean $\underbrace{A\oslash A\oslash \cdots \oslash A}_{k ~ times}.$

If $A\in \M_m$ and $B\in\M_n,$ then $A\oslash B$ can be considered as a block matrix with $i,j-$th block $a_{ii}\J_n+B,~ i=1,2,\ldots ,m,$ in other words $A\oslash B$ is the matrix obtained from $A$ by replacing $a_{i,j}$ by $a_{ii}\J_n+B$. It can be observed that $a_{i,j}+b_{p,q}$ is the $p,q-$th entry of the $i,j-$th block of $A\oslash B.$

All graphs considered here are finite, undirected, connected and simple. The \textit{distance} between two vertices $u,v\in V(G)$ is denoted by $d_{uv}$ and is defined as the length of a shortest path between $u$ and $v$ in $G.$ The \textit{distance matrix} of $G$ is denoted
by $\D(G)$ and is defined by $\D(G)=(d_{uv})_{u,v\in V(G)}.$ Since $\D(G)$ is a real symmetric matrix, all its eigenvalues are real. For a column vector $x=(x_1,\ldots,x_n)^T\in\mathbb{R}^n,$ we have
\beq x^T\D(G)x=\sum_{1\le i<j\le n}d_{ij}x_ix_j. \eeq

The \textit{Wiener index} $W(G)$ of a graph is the sum of the distances between all unordered pairs of vertices of $G$, in other words $W(G)=\frac{S_{\D(G)}}{2}.$ The $\textit{distance spectral radius}~\rho^{\D}(G)$ of $G$ is the largest eigenvalue  of its distance matrix $\D(G).$ 
The transmission, denoted by $Tr(v)$ of a vertex $v$ is the sum of the distances from $v$ to all other vertices in $G$.

The \textit{Cartesian product} $G_1\Box G_2$ of two graphs $G_1$ and $G_2$ is the graph whose vertex set is the Cartesian product $V(G_1)\times V(G_2)$ and in which two vertices $(u, u')$ and $(v, v')$ are adjacent if and only if either $u=v$ and $u'$ is adjacent to $v'$ in $G_2$, or $u'=v'$ and $u$ is adjacent to $v$ in $G_1.$ Let $Gu*Hv$ denote the graph obtained from two graphs $G$ and $H$ by identifying a vertex $u$ from $G$ with a vertex $v$ from $H.$ 

The article have been organized as follows. In \cref{prelim-cartesian}, we discuss some existing results involving Kronecker product of matrices and Cartesian product of graphs. In \cref{traces-cartesian mixed}, we find trace of various compositions of matrices involving Cartesian product. Again in \cref{identities and applications}, we obtain some identities involving Cartesian product of matrices and find some applications in graph theory.

\section{Preliminaries}\label{prelim-cartesian}
Kronecker product has been extensively studied in the literature. Some of the interesting properties of the Kronecker product are given below.

\bl \cite{zha10} \label{tr-kproduct}
If $A\in \M_m$ and $B\in\M_n,$ then $tr(A\otimes B)=tr(A)\times tr(B).$ \el

\bl \cite{zha10} \label{lem-transpose}
If $A\in \M_m$ and $B\in\M_n,$ then $(A\otimes B)^T=A^T\times B^T.$ \el

\bl \cite{zha10} \label{lem-conjugate}
If $A\in \M_m$ and $B\in\M_n,$ then $(A\otimes B)^*=A^*\times B^*.$ \el

\bl \label{cmul} \cite{zha10} For $A\in \M_m,~B\in\M_n,$ and $a,b\in \mathbb{C},$ $aA\otimes bB=abA\otimes B.$   \el 

\bl \cite{zha10} For matrices $A,B,C$ and $D$ of appropriate sizes
\[(A\otimes B)(C\otimes D)=(AC)\otimes (BD).\]
\el 

\bl \cite{zha10} \label{ksimilar} For any $A\in \M_m$ and $B\in \M_n,$ there exist a permutation matrix $P$ such that
\[P^{-1}(A\otimes B)P=B\otimes A.  \] \el 

For more results on Kronecker product, we refer \cite{gra81}. The Cartesian product of two graphs have been studied by many researchers. Here we are interested in Cartesian product of two matrices because for any two connected graphs $G_1$ and $G_2,$ the distance matrix of $G_1\Box G_2$ equals to the Cartesian product of the distance matrices of $G_1$ and $G_2,$  i.e. $\D(G_1\Box G_2)=\D(G_1)\oslash \D(G_2).$ Zhang and Godsil \cite{zg13} found the distance inertia of the Cartesian product of two graphs.

\bt \cite{zg13}
If $G$ and $H$ are two connected graphs, where $V(G)=\{u_1,\ldots, u_m \}$ and $V(H)=\{v_1, \ldots, v_n\},$ then, the inertia of distance matrix of $G\Box H$ is $(n_+(Gu_m *H u_n), (m-1)(n-1)+n_0(Gu_m* Hu_n), n_-(Gu_m*H u_n)).$
\et

\bc \cite{zg13}
Let $T_1$ and $T_2$ be two trees on $m$ and $n$ vertices, respectively. Then the distance inertia of $T_1\Box T_2$ is $(1,(m-1)(n-1),m+n-2).$
\ec

\section{Trace of Cartesian product} \label{traces-cartesian mixed}

Here we consider different compositions and products involving Cartesian product of matrices and evaluate their trace.

\bl \label{tr-cproduct} If $A\in \M_m$ and $B\in\M_n,$ then $tr(A\oslash B)=n.tr(A)+m.tr(B).$ \el 
\pf We have
\beq \mbox{tr}(A\oslash B)&=\mbox{tr}(A\otimes \J_n+\J_n\otimes B)\\
&=\mbox{tr}(A\otimes \J_n)+\mbox{tr}(\J_m\otimes B)\\
&=\mbox{tr}(A)\times \mbox{tr}(\J_n)+\mbox{tr}(\J_m)\times \mbox{tr}(B) \qquad  [\hbox{using } \cref{tr-kproduct} ] \\
&=n.\mbox{tr}(A)+m.\mbox{tr}(B). \eeq \qed

\bt \label{tracemain} If $A_i\in \M_{n_i}$ and $k_i\in \mathbb{C}$ for $i=1,2,\ldots,n,$ then
\[tr(k_1A_1\oslash k_2A_2\oslash\cdots \oslash k_nA_n )=( \Pi_{i=1}^n n_i ) \sum_{i=1}^n \frac{k_itr(A_i)}{n_i}.\] \et 
\pf We prove the result by induction on $n.$ For $n=1,$ there is nothing to prove. For $n=2,$ the result follows from \cref{tr-cproduct}. Suppose the result holds for $n=\ell \le n-1.$ That is 
\be \label{tncpe} \mbox{tr}(k_1A_1\oslash k_2A_2\oslash\cdots \oslash k_\ell A_\ell )=( \Pi_{i=1}^\ell n_i ) \sum_{i=1}^\ell \frac{k_i\mbox{tr}(A_i)}{n_i}.\ee

Now \beq &\mbox{tr}(k_1A_1\oslash k_2A_2\oslash\cdots \oslash k_\ell A_\ell \oslash k_{\ell+1}A_{\ell+1} )\\
=~ &k_{\ell +1} n_{\ell+1}\mbox{tr}(k_1A_1\oslash k_2A_2\oslash\cdots \oslash k_\ell A_\ell)+(\Pi_{i=1}^\ell k_i n_i) \mbox{tr}(A_{\ell+1}) \qquad [ \hbox{using } \cref{tr-cproduct} ] \\
=~&k_{\ell +1} n_{\ell+1}( \Pi_{i=1}^\ell n_i ) \sum_{i=1}^\ell \frac{k_i\mbox{tr}(A_i)}{n_i}+(\Pi_{i=1}^\ell k_i n_i) \mbox{tr}(A_{\ell+1})    \qquad [\hbox{using } \cref{tncpe} ] \\
=~ &( \Pi_{i=1}^{\ell+1} n_i ) \sum_{i=1}^{\ell+1} \frac{k_i \mbox{tr}(A_i)}{n_i}.
\eeq
Hence the result follows by induction. \qed

As immediate corollary of the above theorem we get the following result.

\bc For $A\in \M_n,~ tr(A^{[k]})=k.n^{k-1}tr(A).$ \ec

\bp If $A,B\in \M_n,$ then 
\[tr\big((A+B)\oslash (A-B)\big)=2n.tr(A).  \]
\ep 
\pf From \cref{tr-cproduct}, we have
\beq \mbox{tr}\big((A+B)\oslash (A-B)\big)&=n.\mbox{tr}(A-B)+n.\mbox{tr}(A+B)\\
&=n[\mbox{tr}(A)-\mbox{tr}(B)+\mbox{tr}(A)+\mbox{tr}(B)  ]\\
&=2n.\mbox{tr}(A).
\eeq 
\qed 

\bp For $A\in \M_m, B_i\in \M_n; i=1,2,\ldots ,k,$ then
\[tr(A\otimes(B_1\oslash B_2\oslash \cdots \oslash B_k))=n^{k-1}tr(A)\sum_{i=1}^k tr(B_i). \] \ep 
\pf We have
\beq &~~\mbox{tr}(A\otimes(B_1\oslash B_2\oslash \cdots \oslash B_k))\\
&=\mbox{tr}(A).\mbox{tr}(B_1\oslash B_2\oslash \cdots \oslash B_k) \qquad[\mbox{using }\cref{tr-kproduct}]\\
&=\mbox{tr}(A).n^{k-1}\sum_{i=1}^k \mbox{tr}(B_i) \qquad[\mbox{using }\cref{tracemain}] \\
&=n^{k-1}\mbox{tr}(A)\sum_{i=1}^k \mbox{tr}(B_i).
\eeq 
\qed

\bt If $A_i\in \M_{n_i}$ for $i=1,2,\ldots,t,$ then
\beq tr\big[(A_1\oslash A_2 \oslash\cdots \oslash A_\ell)\otimes(A_{\ell +1}\oslash A_{\ell +2}\oslash \cdots \oslash A_m )\otimes\cdots \otimes (A_r\oslash A_{r+1}\cdots A_t ) \big]\\
= \Pi_{p=1}^t n_p \sum_{i=1}^\ell \frac{tr(A_i)}{n_i}\sum_{j=\ell+1}^m \frac{tr(A_j)}{n_j}\cdots \sum_{k=r}^t \frac{tr(A_k)}{n_k}.   \eeq
\et 
\pf By repeated application of \cref{tr-kproduct} we get 
\beq 
 &~~~\mbox{tr}\big[(A_1\oslash A_2 \oslash\cdots \oslash A_\ell)\otimes(A_{\ell +1}\oslash A_{\ell +2}\oslash \cdots \oslash A_m )\otimes\cdots \otimes (A_r\oslash A_{r+1}\cdots A_t ) \big]\\
&=\big[\mbox{tr}(A_1\oslash A_2\cdots \oslash A_\ell) \big] \big[\mbox{tr}(A_{\ell+1}\oslash A_{\ell+2}\oslash \cdots \oslash A_m ) \big]\cdots \big[\mbox{tr}(A_r\oslash A_{r+1}\oslash \cdots \oslash A_t )  \big]\\
&= \Big(\Pi_{i=1}^\ell \sum_{i=1}^\ell \frac{\mbox{tr}(A_i)}{n_i} \Big)\Big(\Pi_{\ell+1}^m n_j \sum_{j=\ell+1}^m \frac{\mbox{tr}(A_j)}{n_j} \Big)\cdots \Big(\Pi_{k=r}^t n_k \sum_{k=r}^t \frac{\mbox{tr}(A_k)}{n_k} \Big) \quad[ \mbox{using } \cref{tracemain}]  \\
&=\Pi_{p=1}^t n_p \sum_{i=1}^\ell \frac{\mbox{tr}(A_i)}{n_i}\sum_{j=\ell+1}^m \frac{\mbox{tr}(A_j)}{n_j}\cdots \sum_{k=r}^t \frac{\mbox{tr}(A_k)}{n_k}.
\eeq 
\qed 

\bt If $A_i\in \M_{n_i};i=1,2,\ldots ,t,$ then
\beq tr\big[(A_1\otimes A_2 \otimes \cdots \otimes A_\ell)\oslash (A_{\ell+1}\otimes A_{\ell +2}\otimes \cdots \otimes A_m  )\oslash \cdots \oslash (A_r\otimes A_{r+1}\otimes \cdots \otimes A_t)  \big]\\
=\Pi_{p=1}^t n_p \Big[\Pi_{i=1}^\ell \frac{tr(A_i)}{n_i}+\Pi_{j=\ell+1}^m\frac{tr(A_i)}{n_j}+\cdots +\Pi_{k=r}^t \frac{tr(A_k)}{n_k} \Big]. \eeq
\et 
\pf Using \cref{tracemain} and then \cref{tr-kproduct}, we get 
\beq &~~~~ \mbox{tr}\big[(A_1\otimes A_2 \otimes \cdots \otimes A_\ell)\oslash (A_{\ell+1}\otimes A_{\ell +2}\otimes \cdots \otimes A_m  )\oslash \cdots \oslash (A_r\otimes A_{r+1}\otimes \cdots \otimes A_t)  \big]\\
&=\Pi_{p=1}^t n_p\Big[ \frac{\mbox{tr}(A_1\otimes A_2\otimes \cdots A_\ell)}{\Pi_{i=1}^\ell n_i }+\frac{\mbox{tr}(A_{\ell+1}\otimes A_{\ell+2} \cdots \otimes A_m)}{\Pi_{j=\ell+1}^mn_j }+\cdots + \frac{\mbox{tr}(A_r\otimes \cdots \otimes A_t)}{\Pi_{k=r}^tn_k } \Big]\\
&= \Pi_{p=1}^tn_p \Big[\frac{\Pi_{i=1}^\ell \mbox{tr}(A_i) } {\Pi_{i=1}^\ell n_i}+\frac{\Pi_{j=\ell+1}^m \mbox{tr}(A_i)}{\Pi_{j=1}^m n_j}+\cdots + \frac{\Pi_{k=r}^t \mbox{tr}(A_k)}{\Pi_{k=r}^t n_k}  \Big]  \\
&=\Pi_{p=1}^t n_p \Big[\Pi_{i=1}^\ell \frac{\mbox{tr}(A_i)}{n_i}+\Pi_{j=\ell+1}^m\frac{\mbox{tr}(A_i)}{n_j}+\cdots +\Pi_{k=r}^t \frac{\mbox{tr}(A_k)}{n_k} \Big]. \eeq
\qed

\section{Some identities and applications}\label{identities and applications}

From the definition of Cartesian product of two matrices, we get following remarks.
\br If $A$ and $B$ are square matrices and $k\in \mathbb{C},$ then $kA\oslash kB=k(A\oslash B).$ \er

\br For $A\in \M_n$ and any $k\in \mathbb{C},$ $k\oslash A=A+k\J_n=A\oslash k.$ \er


For any square matrices $A$ and $B,$ from the definitions of Kronecker product and Cartesian product, it can be observed that if $a_{i,j}b_{p,q}$ is an entry of $A\otimes B$ then the corresponding entry of $A\oslash B$ is $a_{i,j}+b_{p,q}.$ Thus from \cref{ksimilar}, we see that if $P^{-1}(A\otimes B)P=B\otimes A,$ then for the same $P,$ we get  $P^{-1}(A\oslash B)P=B\oslash A.$ Thus we get the following result.

\br If $A$ and $B$ are square matrices, then $A\oslash B$ is permutation similar to $B\oslash A.$
\er 

\bp \label{prop-transpose} For $A\in \M_m,B\in \M_n,~(A\oslash B)^T=A^T\oslash B^T.$ \ep
\pf By definition we have \[A\oslash B=A\otimes\J_n+J_m\otimes B \] which implies
\beq
(A\oslash B)^T&= (A\otimes\J_n+J_m\otimes B)^T \\
&= (A\otimes \J_n)^T+(\J_m\otimes B)^T \\
&= A^T\otimes\J_n^T+\J_m^T\otimes B^T  \qquad[ \mbox{using \cref{lem-transpose}} ] \\
&=A^T\otimes \J_n+\J_m \otimes B^T \\
&= A^T \oslash B^T.
\eeq 
Hence the result. \qed

By repeated application of \cref{prop-transpose}, we get the following result as a corollary.

\bc For square matrices $A_i$ for $i=1,2,\ldots, n,$
\[(A_1\oslash A_2\oslash \cdots \oslash A_n)^T=A_1^T\oslash A_2^T \oslash \cdots \oslash A_n^T .\] \ec 

Proceeding as in \cref{prop-transpose} and using \cref{lem-conjugate}, we get the following result.

\bp \label{prop-conjugate} For $A\in \M_m,~B\in \M_n,~(A\oslash B)^{*}=A^*\oslash B^*.$
\ep 

By repeated application of \cref{prop-conjugate}, we get the following result as a corollary.

\bc For square matrices $A_i$ for $i=1,2,\ldots, n,$
\[(A_1\oslash A_2\oslash \cdots \oslash A_n)^*=A_1^*\oslash A_2^* \oslash \cdots \oslash A_n^* .\] \ec

\bt \label{csymmetry} If $A\in \M_m, B\in \M_n,$ then $A\oslash B$ is symmetric if and only if $A$ and $B$ are both symmetric.
\et 
\pf If $A$ and $B$ are both symmetric, then $A^T=A$ and $B^T=B.$
Now \beq (A\oslash B)^T&= A^T \oslash B^T \qquad [ \hbox{by }\cref{prop-transpose} ] \\
&= A\oslash B. \eeq 
Therefore $A\oslash B$ is symmetric.

Conversely, suppose that $A\oslash B$ is symmetric. Then $1,1$ block of $A\oslash B$ must be symmetric. But $1,1$ block of $A\oslash B$ is $a_{1,1}\J_n+B$ which is symmetric if and only if $B$ is symmetric. Again since $A\oslash B$ is symmetric, the $1,1$ entry of any $i,j-$th block of $A\oslash B$ must be same as $1,1$ entry of $j,i-$th block of $A\oslash B.$ That is $a_{i,j}+b_{1,1}=a_{j,i}+b_{1,1} \hbox{ for all } i,j=1,2\ldots ,n.$ Which implies that A is symmetric. \qed 

\bt If $A\in \M_m, B\in \M_n,$ then $A\oslash B$ is skew-symmetric if and only if $A$ and $B$ are both skew-symmetric.
\et 
\pf If $A$ and $B$ are both skew-symmetric, then $A^T=-A$ and $B^T=-B.$
Now \beq (A\oslash B)^T &=A^T\oslash B^T  \qquad [ \hbox{by }\cref{prop-transpose} ] \\
&=(-A)\oslash (-B) \\
&= (-A)\otimes \J_n+\J_m\otimes(-B)\\
&=-A\otimes \J_n-\J_n\otimes B \qquad[\hbox{by } \cref{cmul} ] \\
&=-A\oslash B. \eeq 
Therefore $A\oslash B$ is skew-symmetric.

The other direction is similar to that of the proof of \cref{csymmetry}.

\bt If $A\in \M_m$ and $B\in \M_n,$ the $A\oslash B$ is a diagonal matrix if and only if $A=k\J_m$ and $B=-k\J_n$ for some $k\in \mathbb{C}.$ Furthermore in that case $A\oslash B= \large 0.$
\et 
\pf If $A=k\J_m$ and $B=-k\J_n$ for some $k\in \mathbb{C},$ then
\beq A\oslash B &= k\J_m\otimes \J_n +\J_m \otimes (-k\J_n)\\
&= \large 0. \qquad[ \mbox{using } \cref{cmul} ]
\eeq 

Again if $A\oslash B$ is a diagonal matrix, then we must have
\beq 
&a_{i,i}+b_{p,q}=0 \mbox{ for } i=1,2,\ldots ,m \mbox{ and } p,q=1,2,\ldots ,n;~ p\ne q,\\
&a_{i,j}+b_{p,p}=0 \mbox{ for } i,j=1,2,\ldots ,m; ~ i\ne j \mbox{ and } p,q=1,2,\ldots ,n; \\
&a_{i,j}+b_{p,q}=0 \mbox{ for } i=1,2,\ldots ,m; ~i\ne j \mbox{ and } p=1,2,\ldots ,n; ~p\ne q.
\eeq
Solving all those equations we see that all entries of $A$ are equal (say $k$) and all entries of $B$ are also equal ($-k$). Thus we get our required result.
\qed 

\bc There exist no square matrices $A,B$ such that $A\oslash B=I.$ \ec 

\bt If $A,C\in \M_m,~ B,D\in \M_n,$ then $A\oslash B=C\oslash D$ if and only if $C=A-k\J_m$ and $D=B+k\J_n$ for some $k\in \mathbb{C}.$
\et 
\pf If $C=A-k\J_m$ and $D=B+k\J_n$ for some $k\in \mathbb{C},$ then 
\beq 
C\oslash D&= (A-k\J_m)\oslash (B+k\J_n)\\
& = (A-k\J_m)\otimes \J_n +\J_m \otimes (B+k\J_n) \\
& =A\otimes \J_m -k\J_m \otimes \J_n +\J_m \otimes B+\J_m \otimes k\J_n \\
&=A\oslash B.
\eeq
 Conversely, suppose that $A\oslash B= C\oslash D.$ Then every block of $A\oslash B$ equals to the corresponding block of $C\oslash D,$
 \[ \mbox{i.e.  } a_{i,j}\J_n+B=c_{i,j}\J_n+D  \mbox{ for } i,j=1,2,\ldots ,m.\]
 Which implies that $a_{i,j}+b_{p,q}=c_{i,j}+d_{p,q}$ for any $i,j=1,2.\ldots ,m$ and $p,q=1,2,\ldots , n.$
 That is $a_{i,j}-c_{i,j}=d_{p,q}-b_{p,q}$ for any $i,j=1,2.\ldots ,m$ and $p,q=1,2,\ldots , n.$
 Therefore we must have $A-C=\lambda \J_m$ and $D-B=\lambda \J_n$ for some $\lambda \in \mathbb{C}.$
 Hence the theorem follows.
\qed

\bt If $A,B\in \M_n,$ then $A\oslash B=B\oslash A$ if and only if $B=A+k\J_n$ for some $k\in \mathbb{C}.$\et

\pf If $B=A+k\J_n,$ then by direct calculation we have
\[A\oslash B=B\oslash A=A\oslash A+k\J_{n^2}.  \]

Now suppose $A\oslash B=B\oslash A.$ Then
$a_{i,j}+b_{p,q}=b_{i,j}+a_{p,q} \mbox{ for all } i,j,p,q=1,2\ldots,n. $
Therefore 
\beq 
{\ds \sum_{p,q=1}^n }(a_{i,j}+b_{p,q} )&={\ds \sum_{p,q=1}^n }(b_{i,j}+a_{p,q} )\\
\mbox{ which gives  } n^2a_{i,j}+S_B&= n^2b_{i,j}+S_A\\
\mbox{i.e. } b_{i,j} &= a_{i,j}+\frac{S_B-S_A}{n^2} \hs \mbox{ for all } i,j=1,2,\ldots,n.\eeq 
Thus $B=A+k\J_n$ for $k=\frac{S_B-S_A}{n^2}.$ \qed 

\bt If $A,B,C,D\in \M_n,$ then 
\ben[(i)]
\item $(A\oslash B)(C\oslash D)=AC\oslash BD +A\J_n\otimes \J_nB+\J_nC\otimes D\J_n.$
\item $(A\oslash B)\circ(C\oslash D)=(A\circ C)\oslash (B\circ D)+A\otimes D+C\otimes B.$
\een{}
\et 
\pf (i) We have 
\beq (A\oslash B)(C\oslash D)&=(A\otimes \J_n+\J_n\otimes B)(C\otimes \J_n+\J_n\otimes D)\\
&=(A\otimes \J_n)(C\otimes \J_n)+(A\otimes \J_n)(\J_n \otimes B)+(\J_n\otimes B)(C\otimes \J_n)\\&~~~+(\J_n\otimes B)(\J_n \otimes D)\\
&=AC\otimes\J_{n^2}+A\J_n\otimes\J_n B+\J_nC\otimes B\J_n+\J_{n^2}\otimes BD\\
&=AC\oslash BD +A\J_n\otimes \J_nB+\J_nC\otimes D\J_n.
\eeq 

(ii) Here
\beq (A\oslash B)\circ(C\oslash D)&=(A\otimes \J_n+\J_n \otimes B)\circ(C\otimes \J_n +\J_n \otimes D)\\
&=(A\otimes \J_n)\circ (C\otimes \J_n)+(A\otimes \J_n)\circ (\J_n \otimes D)+(\J_n\otimes B)\circ (C\otimes \J_n)\\
&~~~ +(\J_n \otimes B)\circ (\J_n\otimes D)\\
&=(A\circ C)\otimes \J_n +A\otimes D +C\otimes B + \J_n\otimes (B\circ D)\\
&=(A\circ C)\oslash (B\circ D)+A\otimes D+C\otimes B.
\eeq \qed

\bp For matrices $A,B,C$ of suitable orders,
\[(A+B)\oslash C=\frac{1}{2}[A\oslash C+B\oslash C+(A+B)\otimes \J_n]  \] and
\[A\oslash (B+C)=\frac{1}{2}[A\oslash B+A\oslash C+\J\otimes(B+C)].  \]
\ep 
\pf We prove only the first result as the second one can be proved similarly. If $A$ and $B$ are matrices of same order (say $m$) and matrix $C$ is of order $n,$ then
\be \label{prop-distribution}  (A+B)\oslash C &= (A+B)\otimes \J_n +\J_m\otimes C \nonumber \\
&= A\otimes \J_n+B\otimes \J_n+\J_m\otimes C. \qquad[\mbox{since $\otimes$ is distributive}]
\ee 

From \cref{prop-distribution}, we have
\be\label{prodise1} (A+B)\oslash C= A\otimes \J_n+B\oslash C \ee
and
\be\label{prodise2} (A+B)\oslash C=A\oslash C+B\otimes \J_n.  \ee 

Now adding \cref{prodise1} and \cref{prodise2}, we get
\[2((A+B)\oslash C)=A\oslash C+B\oslash C+(A+B)\otimes \J_n.  \]
Hence the result follows.
\qed

\bt \label{sumcartsum} If $A_i\in \M_m, ~B_i\in \M_n$ for $i=1,2,\ldots,k,$ then
\[\Big(\sum_{i=1}^k A_i  \Big)\oslash \Big( \sum_{i=1}^k B_i  \Big)=\sum_{i=1}^k(A_i\oslash B_i).  \]  \et 

\pf We prove the result by induction on $k.$ For $k=1,$ the result is trivial. For $k=2,$ we have
\be \label{longdistribution}
(A_1+A_2)\oslash (B_1+B_2)&=(A_1+A_2)\otimes \J_n+\J_m\otimes (B_1+B_2) \nonumber \\
&= A_1\otimes \J_n +A_2\otimes \J_n+\J_m\otimes B_1 +\J_m \otimes B_2 \quad[\mbox{since $\otimes$ is distributive}] \nonumber \\
&=A_1\oslash B_1+A_2\oslash B_2.
\ee 
Thus the result holds for $k=2.$ Suppose the identity holds for $k=1,2,\ldots,\ell<k,$ then
\beq 
\Big(\sum_{i=1}^{\ell+1} A_i  \Big)\oslash \Big( \sum_{i=1}^{\ell+1} B_i  \Big)&=\Big(\sum_{i=1}^\ell A_i  \Big)\oslash \Big( \sum_{i=1}^\ell B_i  \Big)+A_{\ell+1}\oslash B_{\ell+1} \hs [\mbox{by }\cref{longdistribution}]\\
&= \sum_{i=1}^\ell(A_i\oslash B_i)+A_{\ell+1}\oslash B_{\ell+1}\hs [\mbox{by induction hypothesis}]\\
&=\sum_{i=1}^{\ell+1}(A_i\oslash B_i). \eeq
Hence the result follows.

Using \cref{sumcartsum} repeatedly, we get the following general result.

\bt For $A_i\in \M_m,B_i\in\M_n,\ldots, C_i\in \M_\ell,$ for $i=1,2\ldots,k,$ then
\[ \Big(\sum_{i=1}^k A_i\Big)\oslash \Big(\sum_{i=1}^kB_i\Big)\oslash \cdots \oslash \Big(\sum_{i=1}^kC_i\Big)=\sum_{i=1}^k(A_i\oslash B_i \oslash\cdots \oslash C_i). \]
\et

\bl \label{sumkroneck} If $A$ and $B$ are any square matrices, then
\[S_{A\otimes B}=S_AS_B. \]
\el 
\pf If $A\in \M_m$ and $B\in \M_n,$ then the $i,j-$th block of $A\otimes B$ is $a_{i,j}B$ and $S_{a_{i,j}B}=a_{i,j}S_B.$
Therefore we get \[ S_{A\oslash B}=S_B \sum_{i,j=1}^m a_{i.j}=S_BS_A. \]
Hence the result follows.\qed 

\bt \label{sumcartesian}  If $A\in \M_m$ and $B\in \M_n,$  then
\[S_{A\oslash B}=n^2S_A+m^2S_B. \]
\et 
\pf We have
\beq S_{A\oslash B}&= S_{A\otimes \J_n +\J_m\otimes B}\\
&= S_{A\otimes \J_n}+S_{\J_m\otimes B}\\
&= S_A\times n^2 +m^2\times S_B.   \qquad[ \mbox{using } \cref{sumkroneck}].
\eeq 
Hence the theorem holds. \qed 

As a corollary of \cref{sumcartesian}, we get the expression for the Wiener index of Cartesian product of two connected graphs.

\bc If $G_1$ and $G_2$ are two connected graphs of order $m$ and $n$ respectively, then
\[W(G_1\Box G_2)=n^2W(G_1)+m^2W(G_2). \]
\ec 

As an application of above corollary we get the the following result.

\bc If $H$ is any fixed connected graph and $G_1, G_2$ are connected graphs of same order with $W(G_1)\ge W(G_2),$ then
\[W(H\Box G_1)\ge W(H\Box G_2), \]
with equality if and only if $W(G_1)=W(G_2).$
\ec

\bt \label{cons-rsum-cartesian} If $A\in \M_m$ and $B\in \M_n,$ then $A\oslash B$ has constant row sum if and only if $A$ and $B$ both have constant row sums. \et 
\pf Let us consider any arbitrary row of $A\oslash B.$ If the first entry of that row is $a_{i,1}+b_{j,i},$ then the row sum of that row of $A\oslash B$ equals to
\be \label{crsce} (na_{i,1}+B_j)+(na_{i,2}+B_j)+\cdots +(na_{i,m}+B_j)=nA_j+mB_j.
\ee

Now if $A$ and $B$ have constant row sums, then $A_i=\frac{S_A}{m}$ and $B_i=\frac{S_B}{n}.$ Therefore, by \cref{crsce}, $A\oslash B$ has constant row sum equal to $\frac{n}{m}S_A+\frac{m}{n}S_B.$

Again if $A\oslash B$ has constant row sum (say $k$), then from \cref{crsce} we get 
\[nA_i+mB_j=k \mbox{ for } i=1,2,\ldots , m  \mbox{ and } j=1,2,\ldots ,n. \]
Keeping $i$ fixed, we see that $B_j$ is constant for $j=1,2,\ldots ,n.$ Similarly, keeping $j$ fixed we get $A_i$ is constant for $i=1,2,\ldots ,m.$ Hence, the theorem holds.
\qed 

The following result is a reformulation of \cref{cons-rsum-cartesian}. Therefore, the proof is omitted.

\bt If $A\in \M_m$ and $B\in \M_n,$ then $\1_{mn}$ is an eigenvector of $A\oslash B$ if and only if $\1_m$ and $\1_n$ are eigenvectors of $A$ and $B$ respectively. 
\et

As an application of \cref{cons-rsum-cartesian}, we get the following result as a corollary.

\bc The Cartesian product $G_1\Box G_2$ of two connected graphs $G_1$ and $G_2$ is transmission regular if and only if $G_1$ and $G_2$ are both transmission regular.
\ec 

From the proof of \cref{cons-rsum-cartesian}, we get a lower bound for the distance spectral radius of the Cartesian product of two connected graphs.

\bc If $G_1$ and $G_2$ are two connected graphs of order $m$ and $n$ respectively, then 
\[ \r^\D(G_1\Box G_2)\ge \frac{n}{m}W(G_1)+\frac{m}{n}W(G_2),\]
with equality if and only if $G_1$ and $G_2$ are both transmission regular.
\ec

\end{document}